\documentclass[12pt]{article}%
\usepackage{amsfonts}
\usepackage{sw20bams}
\usepackage{amsmath}
\usepackage{amssymb}
\usepackage{graphicx}%
\providecommand{\U}[1]{\protect\rule{.1in}{.1in}}
\begin{document}

\title{The Maximum of an Asymmetric Simple Random Walk with Reflection}
\author{Steven Finch}
\date{August 24, 2018}
\maketitle

\begin{abstract}
Consider the extreme value of a Bernoulli random walk on the one-dimensional
integer lattice, with reflection at 0, over a finite discrete time interval.
\ Only the asymmetric (biased) case is discussed. \ Asymptotic mean/variance
results are given as the time interval length approaches infinity. \ We
similarly solve an elementary traffic light problem from queueing theory.

\end{abstract}

\footnotetext{Copyright \copyright \ 2018 by Steven R. Finch. All rights
reserved.}Let $X_{0}=0$ and $X_{1}$, $X_{2}$, \ldots, $X_{n}$ be a sequence of
independent random variables satisfying%
\[%
\begin{array}
[c]{ccccccc}%
\mathbb{P}\left(  X_{i}=1\right)  =p, &  & \mathbb{P}\left(  X_{i}=-1\right)
=q, &  & p+q=1, &  & p\leq q
\end{array}
\]
for each $1\leq i\leq n$. \ Define $S_{0}=X_{0}$ and%
\[
S_{j}=\left\{
\begin{array}
[c]{lll}%
\left\vert S_{j-1}+X_{j}\right\vert  &  & \text{strong reflection at the
origin,}\\
\max\left\{  S_{j-1}+X_{j},0\right\}  &  & \text{weak reflection at the
origin}%
\end{array}
\right.
\]
for each $1\leq j\leq n$. \ The simple reflected random walk $S_{0}$, $S_{1}$,
$S_{2}$, \ldots, $S_{n}$ is symmetric if $p=q$ and asymmetric if $p<q$. \ Let%
\[
M_{n}=\max\limits_{0\leq j\leq n}S_{j}
\]
denote the maximum value of the walk over the time interval $[0,n]$. \ We
shall focus entirely on the asymmetric case; a survey of related results
(including those for a symmetric walk with reflection) appears elsewhere
\cite{Fi1-mxasr}.

For the strong scenario, we have%
\[
\mathbb{E}\left(  M_{n}\right)  \sim\frac{\ln(n)}{\ln(\frac{q}{p})}%
+\frac{\gamma+\ln\left(  \frac{(1-2p)^{2}}{2q^{2}}\right)  }{\ln(\frac{q}{p}%
)}+\frac{1}{2}+\varphi(n)
\]
as $n\rightarrow\infty$ and, for the weak scenario,%
\[
\mathbb{E}\left(  M_{n}\right)  \sim\frac{\ln(n)}{\ln(\frac{q}{p})}%
+\frac{\gamma+\ln\left(  \frac{p(1-2p)^{2}}{q^{2}}\right)  }{\ln(\frac{q}{p}%
)}+\frac{1}{2}+\psi(n).
\]
The symbol $\gamma$ denotes Euler's constant \cite{Fi0-mxasr} ; $\varphi$ and
$\psi$ are periodic functions of $\log_{q/p}(n)$ with period $1$ and of small
amplitude. \ It is not surprising that the strong means are slightly larger
than the weak means (because a weakly reflected walk\ can dwell at the origin
indefinitely). \ This is provably true since $1/2>p$. \ For both scenarios, we
have
\[
\mathbb{V}\left(  M_{n}\right)  \sim\frac{\pi^{2}}{6}\frac{1}{\ln\left(
\frac{q}{p}\right)  ^{2}}+\frac{1}{12}+\omega(n)
\]
as $n\rightarrow\infty$ and the function $\omega$, like $\varphi$ and $\psi$,
is effectively negligible.

Different expressions emerge for a certain traffic light problem
\cite{Fi2-mxasr}. \ Let $\ell\geq1$ be an integer. \ Let $X_{0}=0$ and $X_{1}%
$, $X_{2}$, \ldots, $X_{n}$ be a sequence of independent random variables
satisfying%
\[%
\begin{array}
[c]{ccccc}%
\mathbb{P}\left\{  X_{i}=1\right\}  =p, &  & \mathbb{P}\left\{  X_{i}%
=0\right\}  =q &  & \text{if }i\equiv1,2,\ldots,\ell\,\operatorname{mod}%
\,2\ell;
\end{array}
\]%
\[%
\begin{array}
[c]{ccccc}%
\mathbb{P}\left\{  X_{i}=0\right\}  =p, &  & \mathbb{P}\left\{  X_{i}%
=-1\right\}  =q &  & \text{if }i\equiv\ell+1,\ell+2,\ldots,2\ell
\,\operatorname{mod}\,2\ell
\end{array}
\]
for each $1\leq i\leq n$. \ Define $S_{0}=X_{0}$ and $S_{j}=\max\left\{
S_{j-1}+X_{j},0\right\}  $ for all $1\leq j\leq n$. \ Thus cars arrive at a
one-way intersection according to a Bernoulli($p$) distribution; \ when the
signal is red ($1\leq i\leq\ell$), no cars may leave; \ when the signal is
green ($\ell+1\leq i\leq2\ell$), a car must leave (if there is one). The
quantity $M_{n}=\max\nolimits_{0\leq j\leq n}S_{j}$ is the worst-case traffic
congestion (as opposed to the average-case often cited). \ Only the
circumstance when $\ell=1$ is amenable to rigorous treatment, as far as is
known. \ We have%
\[
\mathbb{E}\left(  M_{n}\right)  \sim\frac{\ln(n)}{2\ln(\frac{q}{p})}%
+\frac{\gamma+\ln\left(  \frac{p(1-2p)^{2}}{2q^{3}}\right)  }{2\ln(\frac{q}%
{p})}+\frac{1}{2}+\psi(n),
\]%
\[
\mathbb{V}\left(  M_{n}\right)  \sim\frac{\pi^{2}}{24}\frac{1}{\ln\left(
\frac{q}{p}\right)  ^{2}}+\frac{1}{12}+\omega(n)
\]
as $n\rightarrow\infty$, assuming $p<q$. \ For $\ell>1$, only computer
simulation-based estimates are available.

Sections 1 and 2 cover generating functions and singularity analyses
corresponding to strongly RRWs and weakly RRWs, respectively. \ Section 5 does
likewise for the traffic light problem ($\ell=1$). \ We focus on the
calculation of moments in Section 3. \ Sections 4 and 6 provide extensive
verification by use of simulation.

\section{Strong Scenario}

For $k=1,2,\ldots$, define $(k+1)\times(k+1)$ matrices%
\[%
\begin{array}
[c]{ccccc}%
A_{1}=\left(
\begin{array}
[c]{cc}%
0 & q\\
1 & 0
\end{array}
\right)  , &  & A_{2}=\left(
\begin{array}
[c]{ccc}%
0 & q & 0\\
1 & 0 & q\\
0 & p & 0
\end{array}
\right)  , &  & A_{3}=\left(
\begin{array}
[c]{cccc}%
0 & q & 0 & 0\\
1 & 0 & q & 0\\
0 & p & 0 & q\\
0 & 0 & p & 0
\end{array}
\right)  ,
\end{array}
\]%
\[%
\begin{array}
[c]{ccccc}%
A_{4}=\left(
\begin{array}
[c]{ccccc}%
0 & q & 0 & 0 & 0\\
1 & 0 & q & 0 & 0\\
0 & p & 0 & q & 0\\
0 & 0 & p & 0 & q\\
0 & 0 & 0 & p & 0
\end{array}
\right)  , &  & A_{5}=\left(
\begin{array}
[c]{cccccc}%
0 & q & 0 & 0 & 0 & 0\\
1 & 0 & q & 0 & 0 & 0\\
0 & p & 0 & q & 0 & 0\\
0 & 0 & p & 0 & q & 0\\
0 & 0 & 0 & p & 0 & q\\
0 & 0 & 0 & 0 & p & 0
\end{array}
\right)  , &  & \ldots
\end{array}
\]
and column vectors $\varepsilon_{k+1}=(1,0,0,\ldots,0)$, $1_{k+1}%
=(1,1,\ldots,1)$. It is not difficult (e.g., starting with \cite{Fi1-mxasr})
to show that%
\[
\mathbb{P}\left\{  M_{n}\leq k\right\}  =1_{k+1\,}^{^{\prime}}A_{k}%
^{n}\;\varepsilon_{k+1}%
\]
and $A^{n}$ denotes $n^{\text{th}}$ matrix power. \ This yields a generating
function%
\begin{align*}
G_{k}(z) &  =%
{\displaystyle\sum\limits_{n=0}^{\infty}}
\mathbb{P}\left\{  M_{n}\leq k\right\}  z^{n}=1_{k+1\,}^{^{\prime}}%
{\displaystyle\sum\limits_{n=0}^{\infty}}
(A_{k}\,z)^{n}\varepsilon_{k+1}=1_{k+1\,}^{^{\prime}}\left(  I-A_{k}%
\,z\right)  ^{-1}\varepsilon_{k+1}\\
&  =\frac{\det\left(  B_{k}\right)  }{\det\left(  I-A_{k}\,z\right)  }%
=\frac{Q_{k}}{R_{k}}%
\end{align*}
where $B_{k}$ is obtained from $I-A_{k}\,z$ via replacing the first row by
$1_{k+1\,}^{^{\prime}}$. \ Expanding the determinant with respect to the last
row, we find linear recursive formulas
\[%
\begin{array}
[c]{ccccc}%
R_{k}=R_{k-1}-p\,q\,z^{2}R_{k-2}, &  & R_{0}=1, &  & R_{1}=1-q\,z^{2};
\end{array}
\]%
\[%
\begin{array}
[c]{ccccc}%
Q_{k}=Q_{k-1}-p\,q\,z^{2}Q_{k-2}+p^{k-1}z^{k}, &  & Q_{-1}=0, &  & Q_{0}=1.
\end{array}
\]
Setting $t=\sqrt{1-4pqz^{2}}$, explicit solutions are as follows:%
\[
R_{k}=\frac{1}{2t}\left[  u\left(  \frac{1+t}{2}\right)  ^{k}+v\left(
\frac{1-t}{2}\right)  ^{k}\right]  ,
\]%
\[
Q_{k}=\frac{1}{2t(1-z)}\left[  u\left(  \frac{1+t}{2}\right)  ^{k}+v\left(
\frac{1-t}{2}\right)  ^{k}-2tz(pz)^{k}\right]
\]
where%
\begin{align*}
u &  =1+t-2qz^{2}=1+t-2q\left(  \frac{1-t^{2}}{4pq}\right)  \\
&  =\left(  \frac{1+t}{2p}\right)  \left[  2p-(1-t)\right]  =\left(
\frac{1+t}{2p}\right)  \left(  -1+2p+t\right)  ,
\end{align*}%
\begin{align*}
v &  =-1+t+2qz^{2}=-1+t+2q\left(  \frac{1-t^{2}}{4pq}\right)  \\
&  =\left(  \frac{1-t}{2p}\right)  \left[  -2p+(1+t)\right]  =\left(
\frac{1-t}{2p}\right)  \left(  1-2p+t\right)  .
\end{align*}
To assess the asymptotics for the coefficients of $G_{k}(z)$, it suffices to
examine the zero $z_{k}$ of its denominator that is closest to the origin.
\ The equation $R_{k}=0$ can be rewritten as%
\[
\left(  \frac{1+t}{1-t}\right)  ^{k}=\frac{v}{-u}=\frac{1-t}{1+t}\cdot
\frac{1-2p+t}{1-2p-t},
\]
that is,%
\[
\left(  \frac{1+t}{1-t}\right)  ^{k+1}=\frac{1-2p+t}{1-2p-t}.
\]
For suitably large $k$, there is exactly one solution $t_{k}$ of the preceding
equation with positive real part; further, $t_{k}$ is real and satisfies
$0<t_{k}<1$. \ The details (involving Rouch\'{e}'s theorem) are omitted. \ It
follows that $z_{k}$ is real; as both $z_{k}$ and $-z_{k}$ are zeroes of
$R_{k}$, we choose $z_{k}$\ to be positive (without loss of generality).
\ Further,%
\[
\frac{z_{k}-1}{(p/q)^{k}}\rightarrow\frac{(1-2p)^{2}}{2q^{2}}%
\]
for the strong scenario as $k\rightarrow\infty$, which implies that
\cite{HW-mxasr, W-mxasr}%
\[
\mathbb{P}\left\{  M_{n}\leq\log_{q/p}(n)+h\right\}  \sim\exp\left[
-\frac{(1-2p)^{2}}{2q^{2}}\left(  \frac{q}{p}\right)  ^{-h}\right]
\]
as $n\rightarrow\infty$. \ Consequences of such a discrete Gumbel
distributional limit will be explored shortly.

\section{Weak Scenario}

The analog of the $A_{k}$ matrix here is%
\[
\left(
\begin{array}
[c]{ccccccc}%
q & q & 0 & 0 & \cdots & 0 & 0\\
p & 0 & q & 0 & \cdots & 0 & 0\\
0 & p & 0 & q & \cdots & 0 & 0\\
0 & 0 & p & 0 & \cdots & 0 & 0\\
\vdots & \vdots & \vdots & \vdots & \ddots & \vdots & \vdots\\
0 & 0 & 0 & 0 & \cdots & 0 & q\\
0 & 0 & 0 & 0 & \cdots & p & 0
\end{array}
\right)
\]
which gives rise as before to recursions%
\[%
\begin{array}
[c]{ccccc}%
R_{k}=R_{k-1}-pqz^{2}R_{k-2}, &  & R_{-1}=1, &  & R_{0}=1-qz;
\end{array}
\]%
\[%
\begin{array}
[c]{ccccc}%
Q_{k}=Q_{k-1}-pqz^{2}Q_{k-2}+p^{k}z^{k}, &  & Q_{-1}=0, &  & Q_{0}=1.
\end{array}
\]
Solving these, we have%
\[
R_{k}=\frac{1}{2t}\left[  u\left(  \frac{1+t}{2}\right)  ^{k}+v\left(
\frac{1-t}{2}\right)  ^{k}\right]  ,
\]%
\[
Q_{k}=\frac{1}{2t(1-z)}\left[  u\left(  \frac{1+t}{2}\right)  ^{k}+v\left(
\frac{1-t}{2}\right)  ^{k}-2t(pz)^{k+1}\right]
\]
where%
\begin{align*}
u &  =1+t-q(1+t)z-2pqz^{2}=1+t-q(1+t)\sqrt{\frac{1-t^{2}}{4pq}}-2pq\left(
\frac{1-t^{2}}{4pq}\right)  \\
&  =\left(  \frac{1+t}{2}\right)  \left[  2-q\sqrt{\frac{1-t^{2}}{pq}%
}-(1-t)\right]  =\left(  \frac{1+t}{2}\right)  \left(  1+t-\sqrt{\frac{q}{p}%
}\sqrt{1-t^{2}}\right)  ,
\end{align*}%
\begin{align*}
v &  =-1+t+q(1-t)z+2pqz^{2}=-1+t+q(1-t)\sqrt{\frac{1-t^{2}}{4pq}}+2pq\left(
\frac{1-t^{2}}{4pq}\right)  \\
&  =\left(  \frac{1-t}{2}\right)  \left[  -2+q\sqrt{\frac{1-t^{2}}{pq}%
}+(1+t)\right]  =\left(  \frac{1-t}{2}\right)  \left(  -1+t+\sqrt{\frac{q}{p}%
}\sqrt{1-t^{2}}\right)  .
\end{align*}
The equation $R_{k}=0$ can be rewritten as%
\[
\left(  \frac{1+t}{1-t}\right)  ^{k}=\frac{-v}{u}=\frac{1-t}{1+t}\cdot
\frac{1-t-\sqrt{\frac{q}{p}}\sqrt{1-t^{2}}}{1+t-\sqrt{\frac{q}{p}}%
\sqrt{1-t^{2}}},
\]
that is,%
\[
\left(  \frac{1+t}{1-t}\right)  ^{k+1}=\frac{1-t-\sqrt{\frac{q}{p}}%
\sqrt{1-t^{2}}}{1+t-\sqrt{\frac{q}{p}}\sqrt{1-t^{2}}}.
\]
By reasoning akin to earlier,
\[
\frac{z_{k}-1}{(p/q)^{k}}\rightarrow\frac{p(1-2p)^{2}}{q^{2}}%
\]
for the weak scenario as $k\rightarrow\infty$, which implies that%
\[
P\left\{  M_{n}\leq\log_{q/p}(n)+h\right\}  \sim\exp\left[  -\frac
{p(1-2p)^{2}}{q^{2}}\left(  \frac{q}{p}\right)  ^{-h}\right]
\]
as $n\rightarrow\infty$.

\section{Mean and Variance}

Fix $c>0$ and $r>1$. \ Forget temporarily the discrete nature of our
distributional limits. \ To evaluate moments associated with a
\textit{continuous} Gumbel CDF%
\[%
\begin{array}
[c]{ccc}%
\exp\left[  -c\,r^{-y}\right]  =\exp\left[  -e^{-\left(  \ln(r)y-\ln
(c)\right)  }\right]  , &  & -\infty<y<\infty
\end{array}
\]
merely set%
\[
\frac{y-\alpha}{\beta}=\ln(r)y-\ln(c)
\]
i.e.,%
\[%
\begin{array}
[c]{ccc}%
\dfrac{1}{\beta}=\ln(r), &  & \alpha=\dfrac{\ln(c)}{\ln(r)}%
\end{array}
\]
i.e., \cite{G-mxasr}%
\[%
\begin{array}
[c]{ccc}%
\mathbb{E}(Y)=\alpha+\beta\,\gamma=\dfrac{\ln(c)+\gamma}{\ln(r)}, &  &
\mathbb{V}(Y)=\dfrac{\pi^{2}}{6}\beta^{2}=\dfrac{\pi^{2}}{6}\dfrac{1}%
{\ln(r)^{2}}.
\end{array}
\]
Return to the discrete domain is achieved by addition of correction terms:
\[
\mathbb{E}\left(  M_{n}-\log_{r}(n)\right)  \sim\frac{\ln\left(  c\right)
+\gamma}{\ln(r)}+\frac{1}{2}+\varphi(n),
\]%
\[
\mathbb{V}\left(  M_{n}\right)  \sim\dfrac{\pi^{2}}{6}\dfrac{1}{\ln(r)^{2}%
}+\frac{1}{12}+\omega(n)
\]
as $n\rightarrow\infty$, via appropriate generalization of \cite{HL-mxasr,
LP-mxasr} (who unnecessarily restrict $r$ to be exactly $2$). \ The
calculation of higher order moments is also possible.

\section{Walk Data}

Let $n=10^{10}$. \ For each $p\in\{1/5,1/3,1/4,3/7,1/8\}$, we\ generated
$40000$ strongly RRWs and produced an empirical histogram for the maximum
$M_{n}$. \ Figures 1--5 contain these histograms (in blue) along with our
theoretical predictions (in red). \ The fit is excellent. \ 

Similarly, we\ generated $40000$ weakly RRWs and produced a histogram for the
maximum $M_{n}$. \ Figures 6--10 contain these histograms along with our
theoretical predictions. \ The fit, again, is excellent.

Each histogram is accompanied by an experimental mean, mean square and
standard deviation, as well as our theoretical values.

\section{Traffic Light}

The analog of the $A_{k}$ matrix here is%
\[
\left(
\begin{array}
[c]{ccccccc}%
1-p^{2} & q^{2} & 0 & 0 & \cdots & 0 & 0\\
p^{2} & 2pq & q^{2} & 0 & \cdots & 0 & 0\\
0 & p^{2} & 2pq & q^{2} & \cdots & 0 & 0\\
0 & 0 & p^{2} & 2pq & \cdots & 0 & 0\\
\vdots & \vdots & \vdots & \vdots & \ddots & \vdots & \vdots\\
0 & 0 & 0 & 0 & \cdots & 2pq & q^{2}\\
0 & 0 & 0 & 0 & \cdots & p^{2} & pq
\end{array}
\right)
\]
but recursions for $R_{k}$ and $Q_{k}$ are more complicated than those for
walks. \ We have%
\[
R_{k}=\tilde{R}_{k}+p\,q\,z\,\tilde{R}_{k-1},
\]%
\[
Q_{k}=\tilde{Q}_{k}+p\,q\,z\,\tilde{Q}_{k-1}%
\]
and, in turn,%
\[%
\begin{array}
[c]{ccccc}%
\tilde{R}_{k}=(1-2pqz)\tilde{R}_{k-1}-p^{2}q^{2}z^{2}\tilde{R}_{k-2}, &  &
\tilde{R}_{-1}=1, &  & \tilde{R}_{0}=1-\left(  1-p^{2}\right)  z;
\end{array}
\]%
\[%
\begin{array}
[c]{ccccc}%
\tilde{Q}_{k}=(1-2pqz)\tilde{Q}_{k-1}-p^{2}q^{2}z^{2}\tilde{Q}_{k-2}%
+p^{2k}z^{k}, &  & \tilde{Q}_{-1}=0, &  & \tilde{Q}_{0}=1.
\end{array}
\]
\newpage%

\begin{figure}[ptb]%
\centering
\includegraphics[
height=2.9992in,
width=3.6685in
]%
{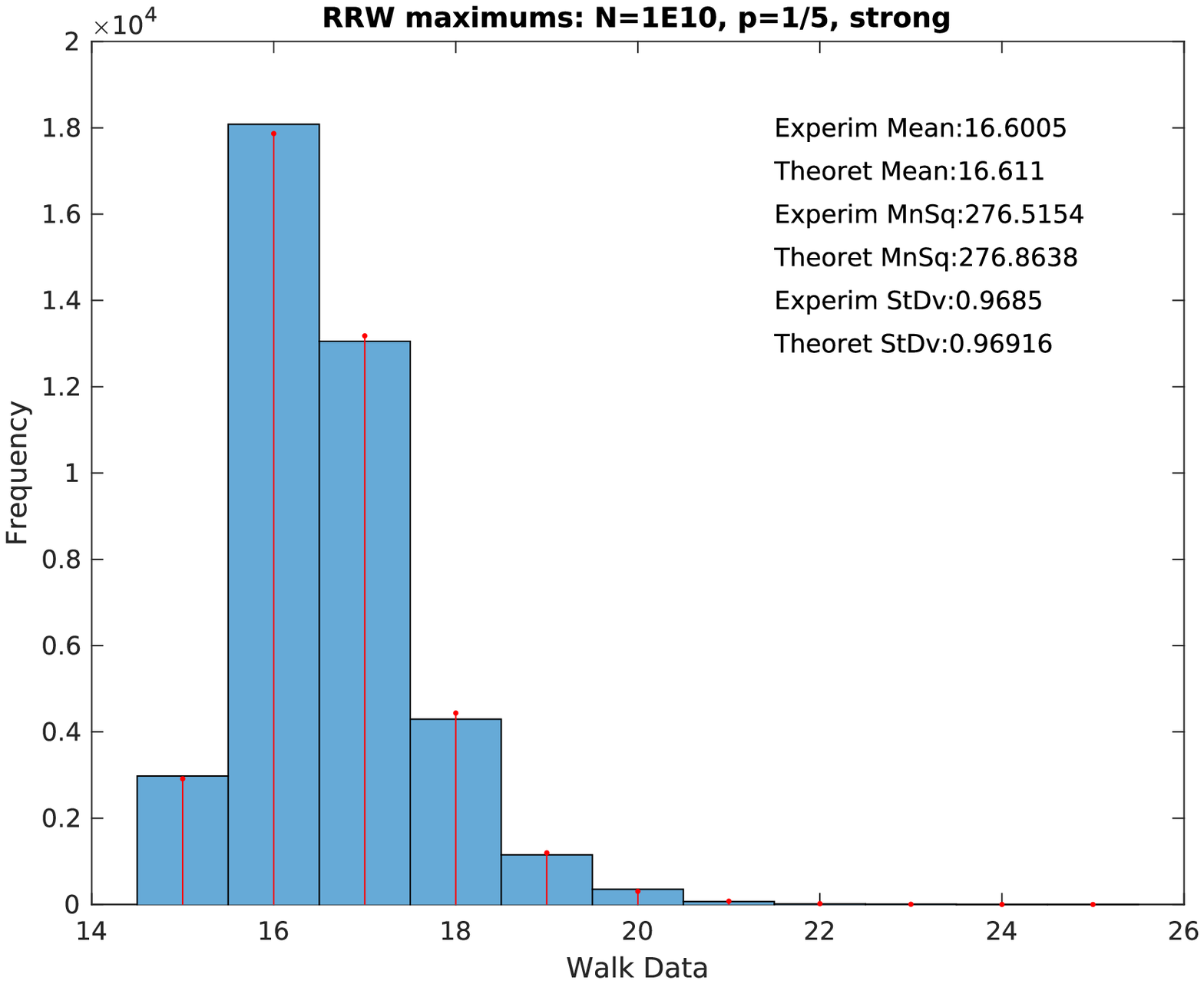}%
\caption{ }%
\end{figure}
\begin{figure}[ptb]%
\centering
\includegraphics[
height=2.9992in,
width=3.7922in
]%
{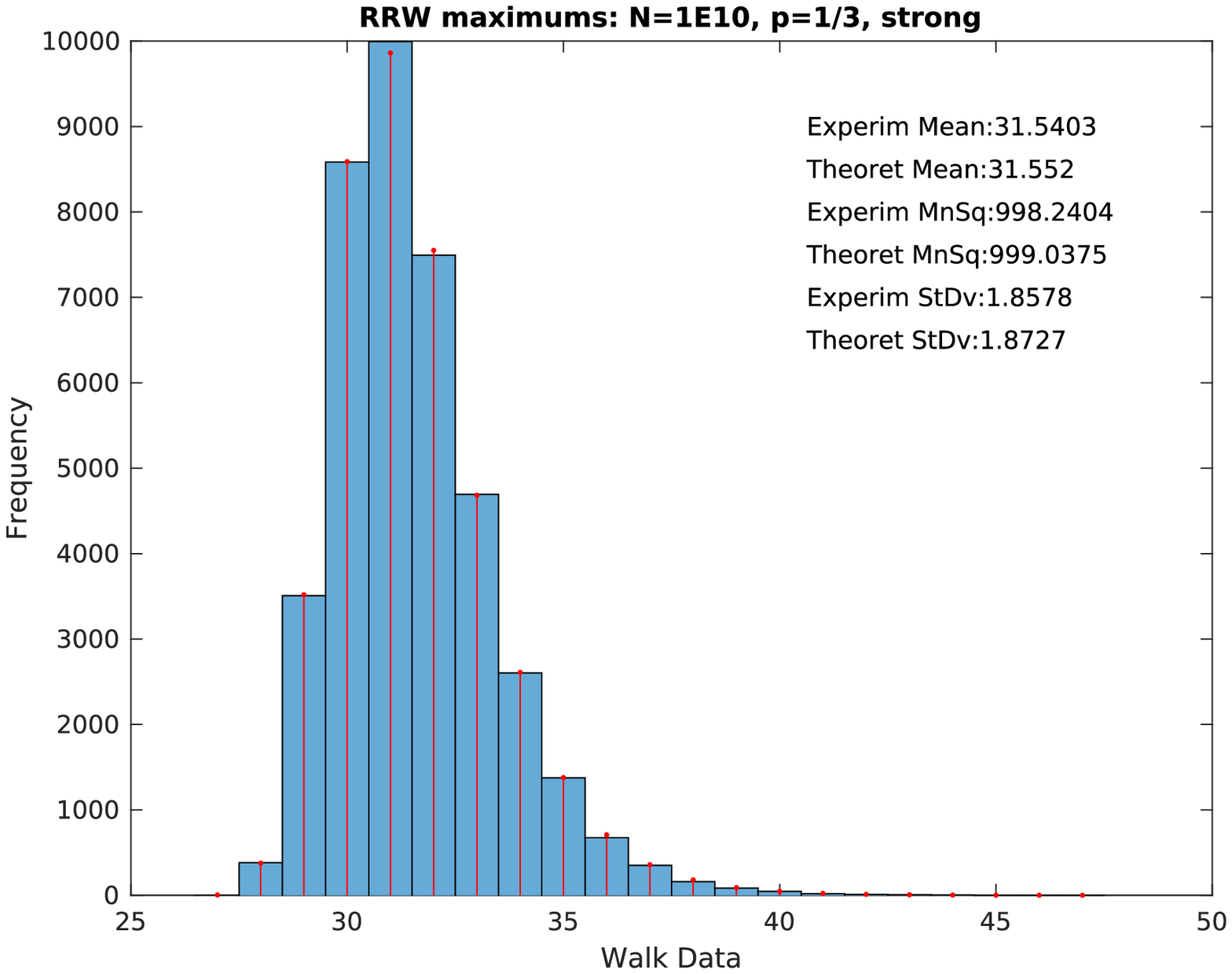}%
\caption{ }%
\end{figure}
\begin{figure}[ptb]%
\centering
\includegraphics[
height=2.9992in,
width=3.7836in
]%
{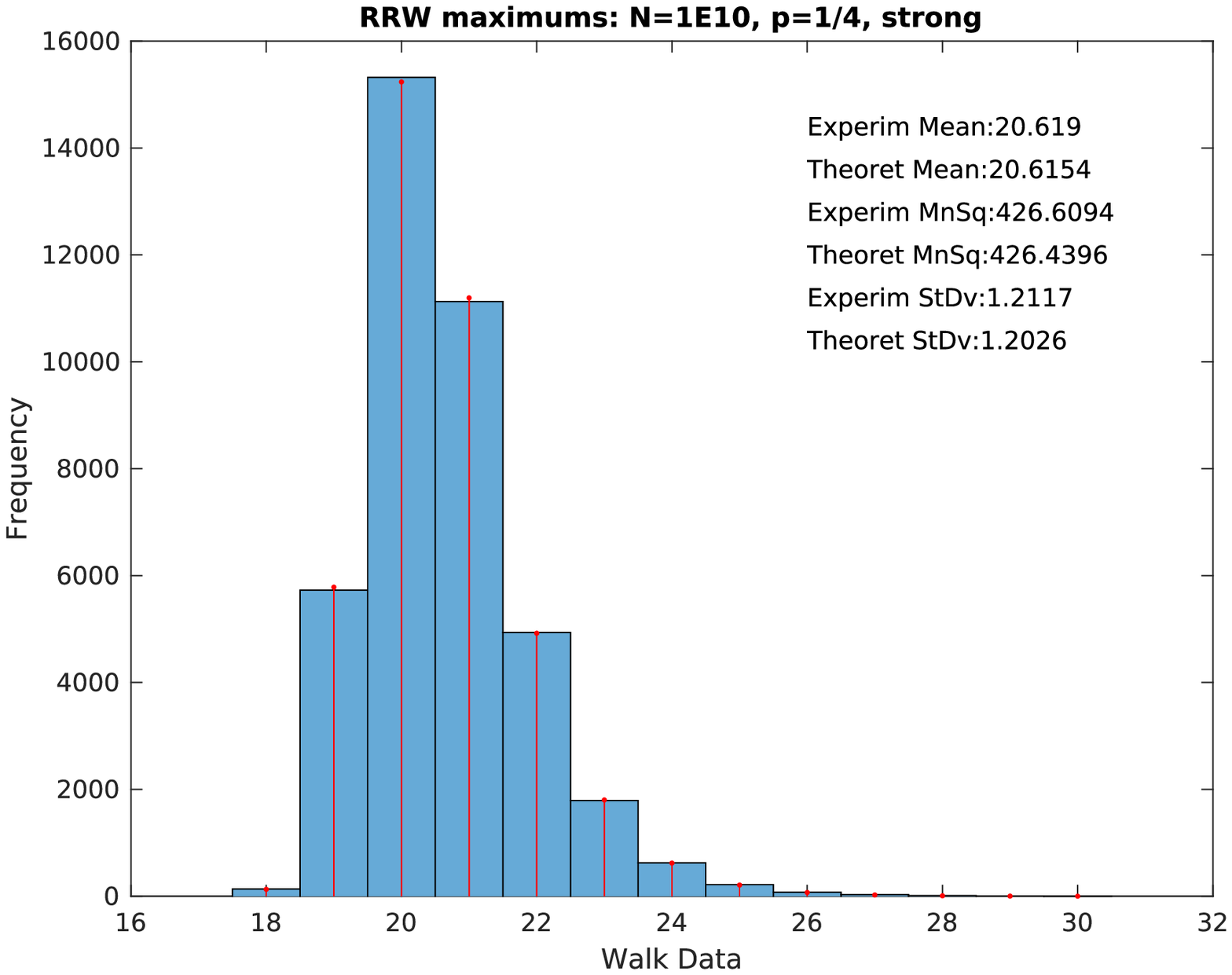}%
\caption{ }%
\end{figure}
\begin{figure}[ptb]%
\centering
\includegraphics[
height=2.9992in,
width=3.7611in
]%
{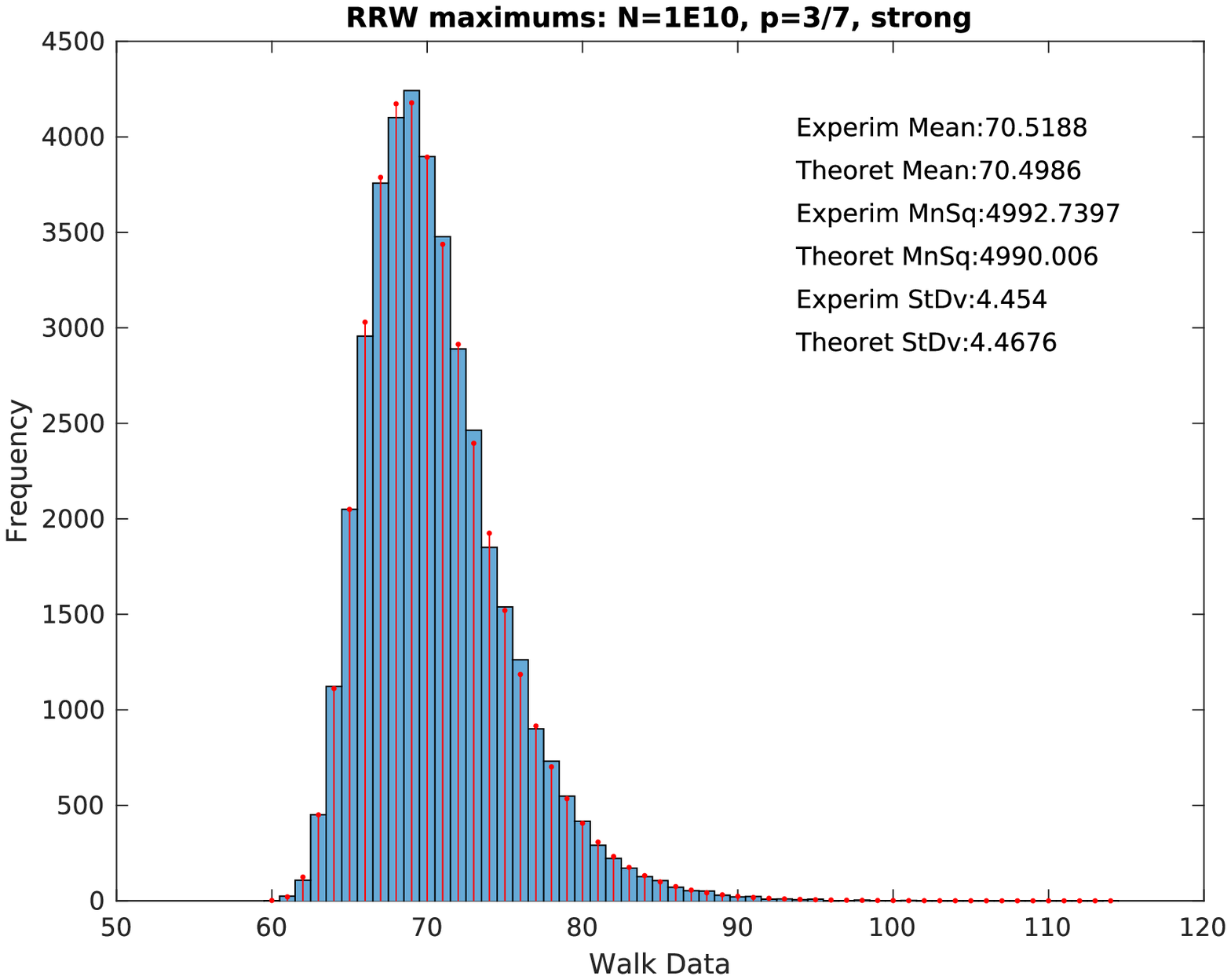}%
\caption{ }%
\end{figure}
\begin{figure}[ptb]%
\centering
\includegraphics[
height=2.9992in,
width=3.6685in
]%
{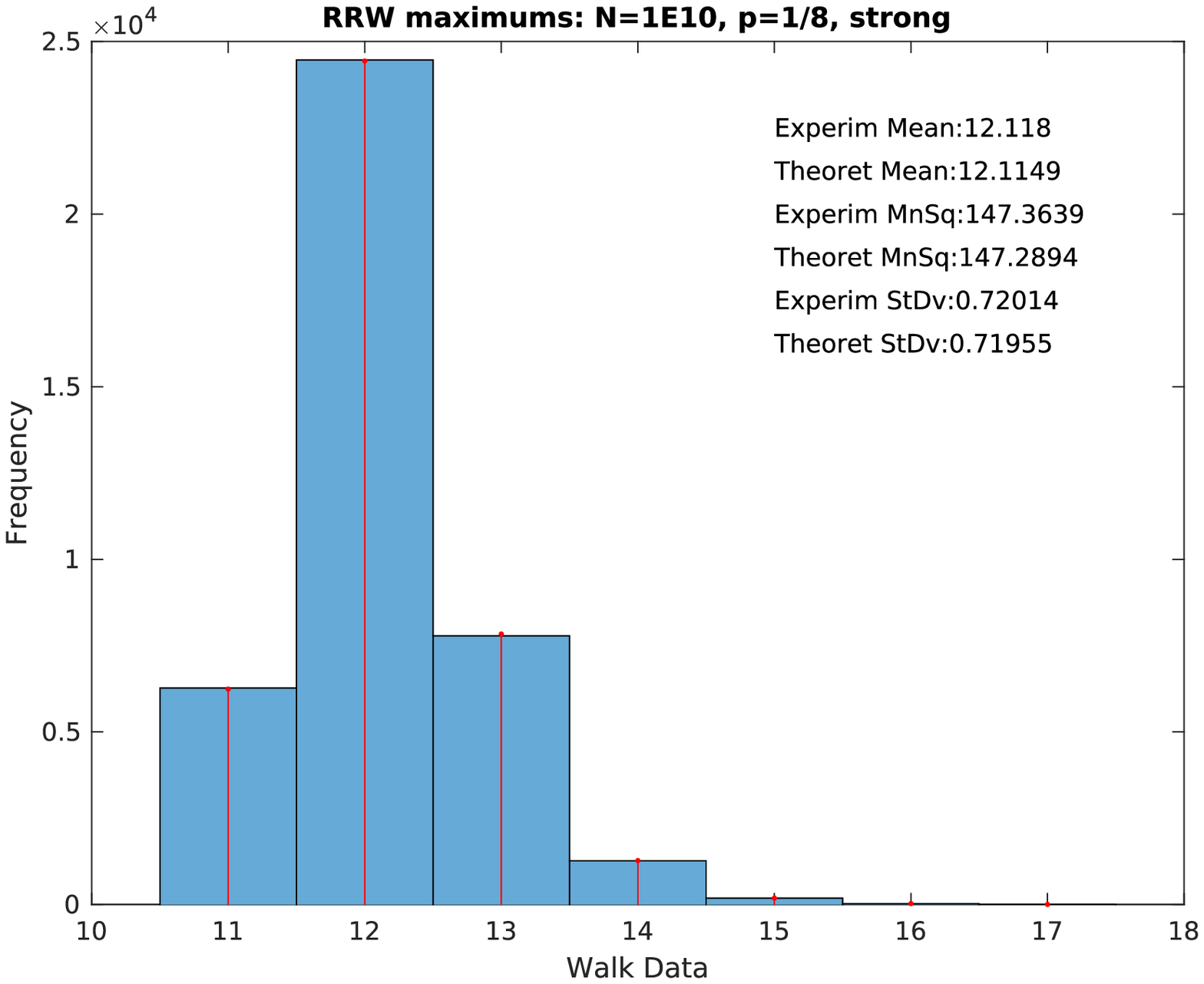}%
\caption{ }%
\end{figure}
\begin{figure}[ptb]%
\centering
\includegraphics[
height=2.9992in,
width=3.7922in
]%
{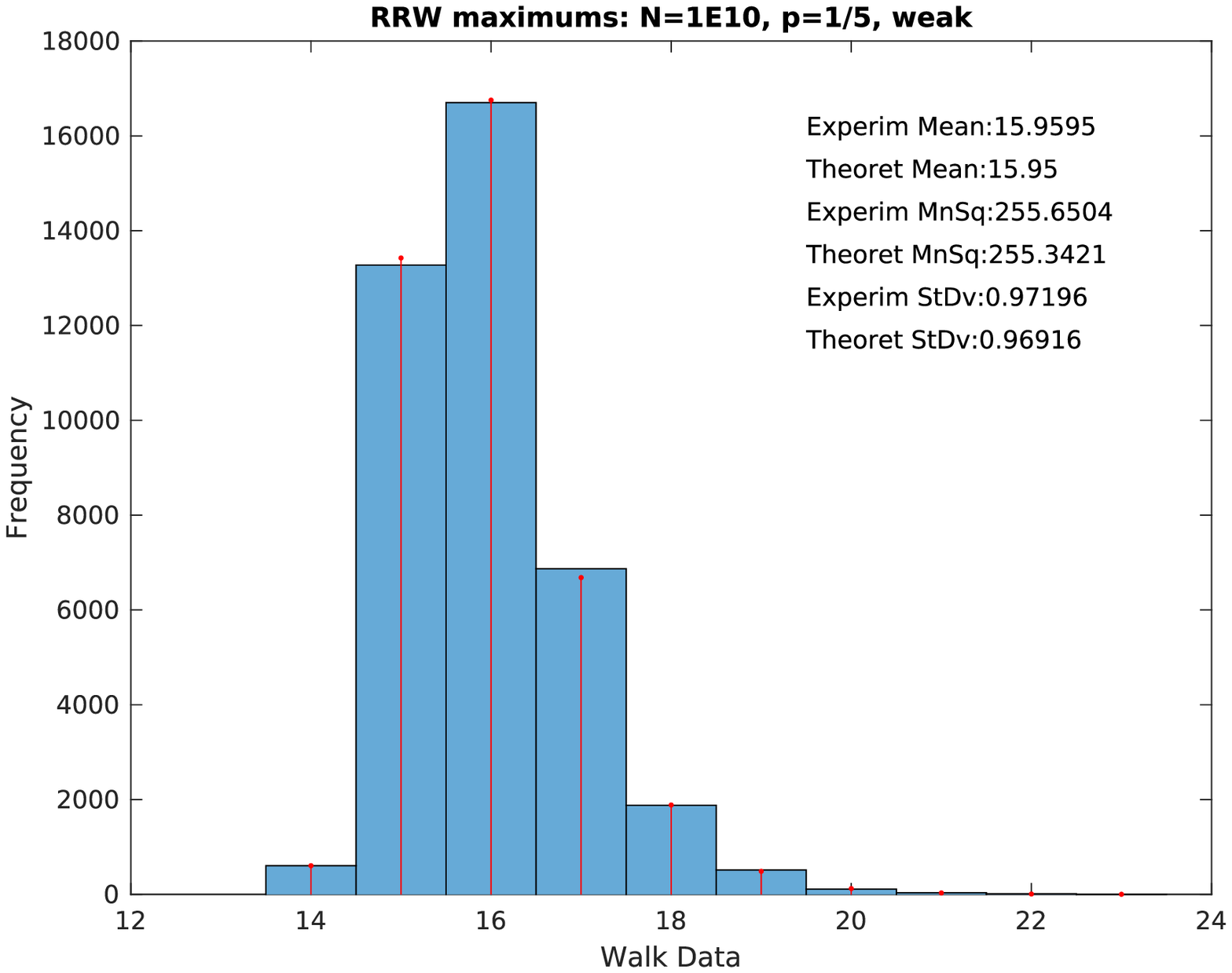}%
\caption{ }%
\end{figure}
\begin{figure}[ptb]%
\centering
\includegraphics[
height=2.9992in,
width=3.7836in
]%
{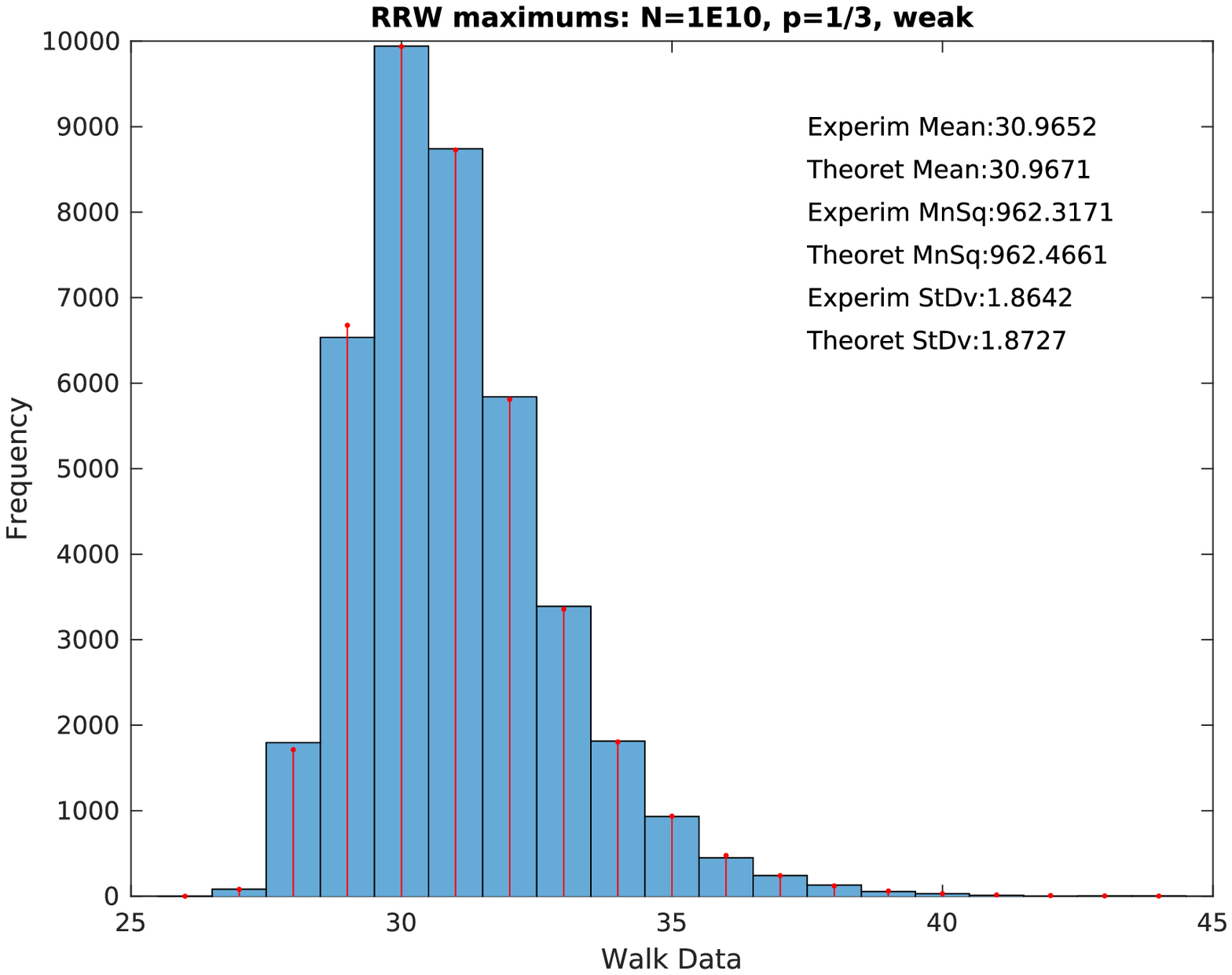}%
\caption{ }%
\end{figure}
\begin{figure}[ptb]%
\centering
\includegraphics[
height=2.9992in,
width=3.7922in
]%
{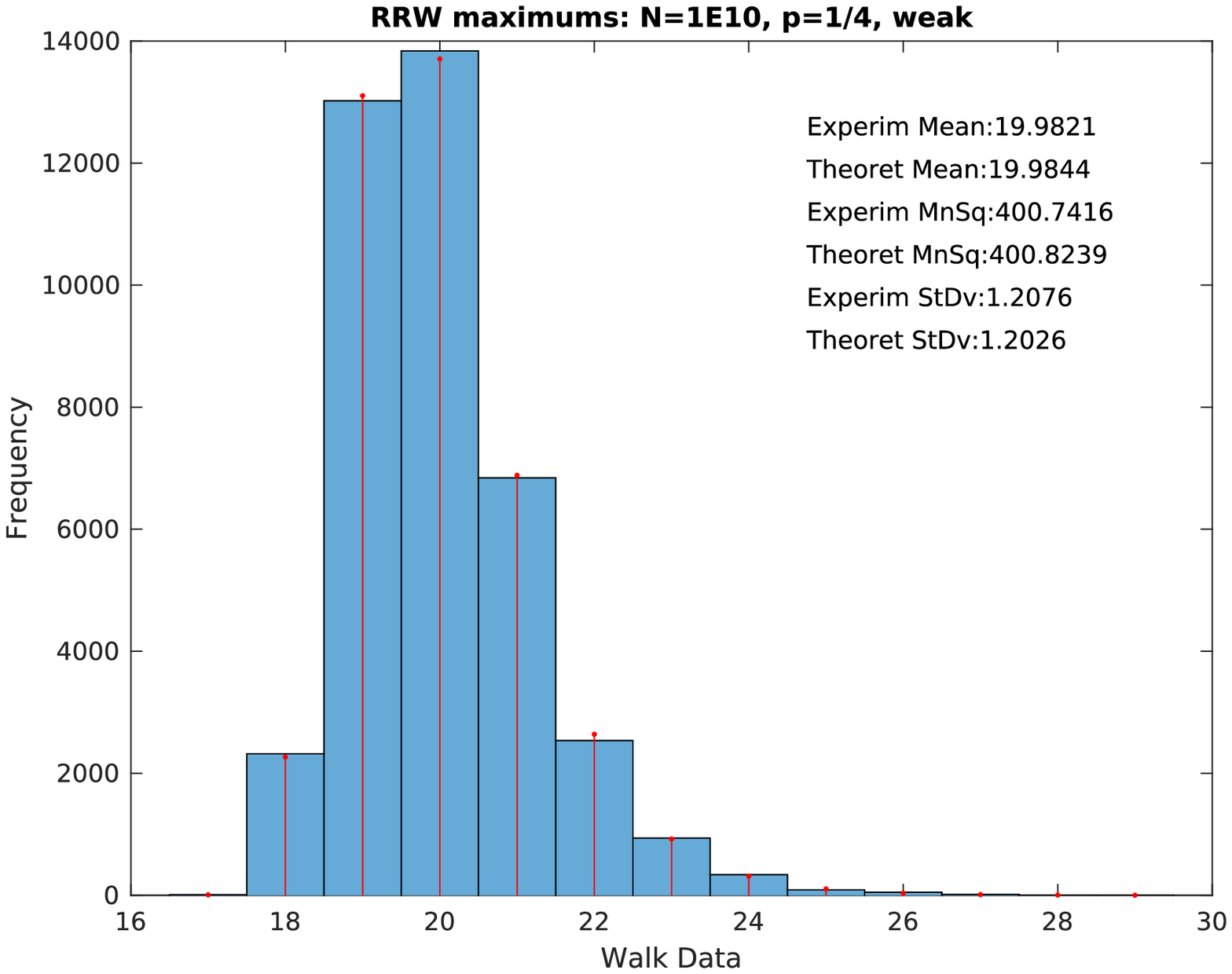}%
\caption{ }%
\end{figure}
\begin{figure}[ptb]%
\centering
\includegraphics[
height=2.9992in,
width=3.7611in
]%
{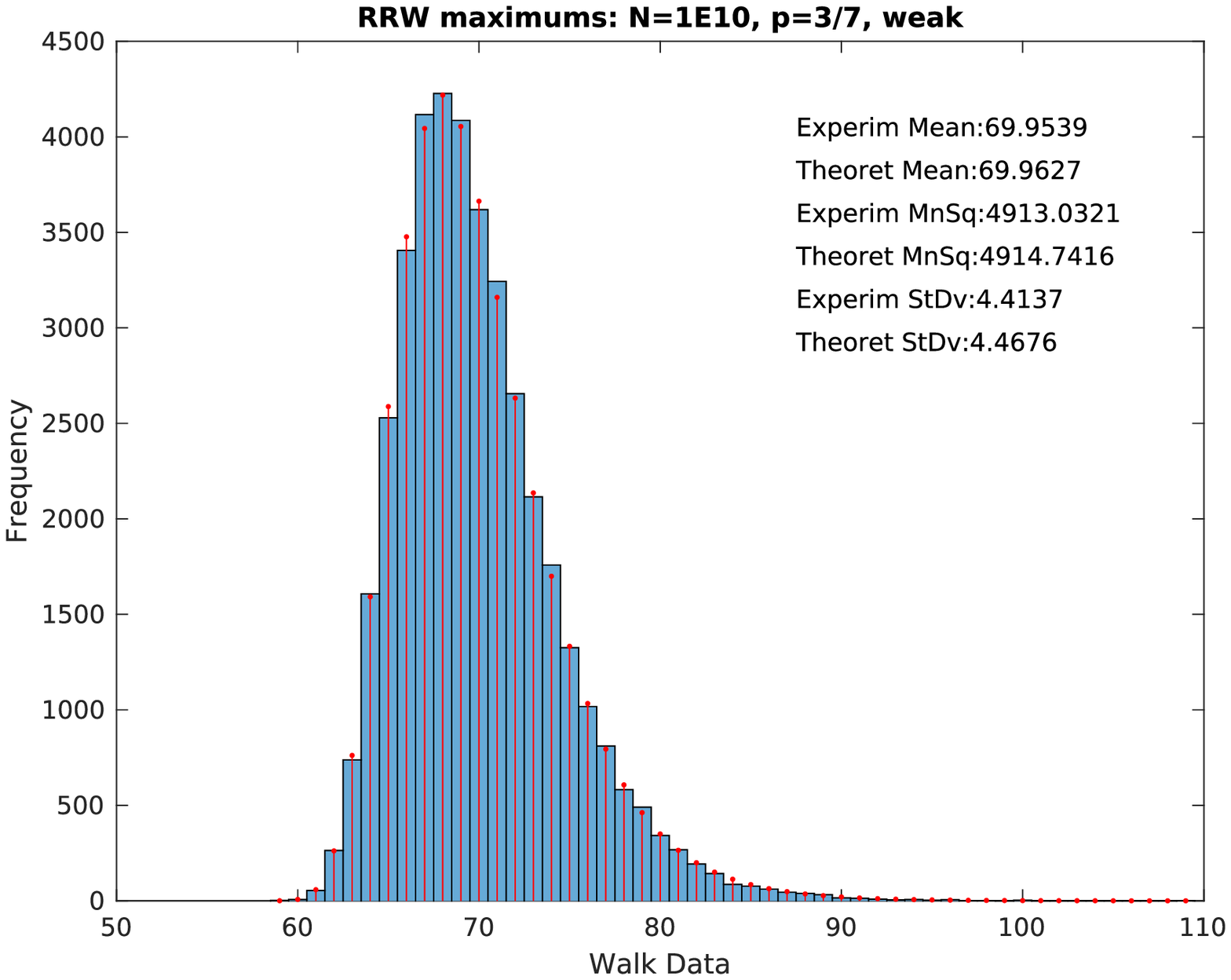}%
\caption{ }%
\end{figure}
\begin{figure}[ptb]%
\centering
\includegraphics[
height=2.9992in,
width=3.6685in
]%
{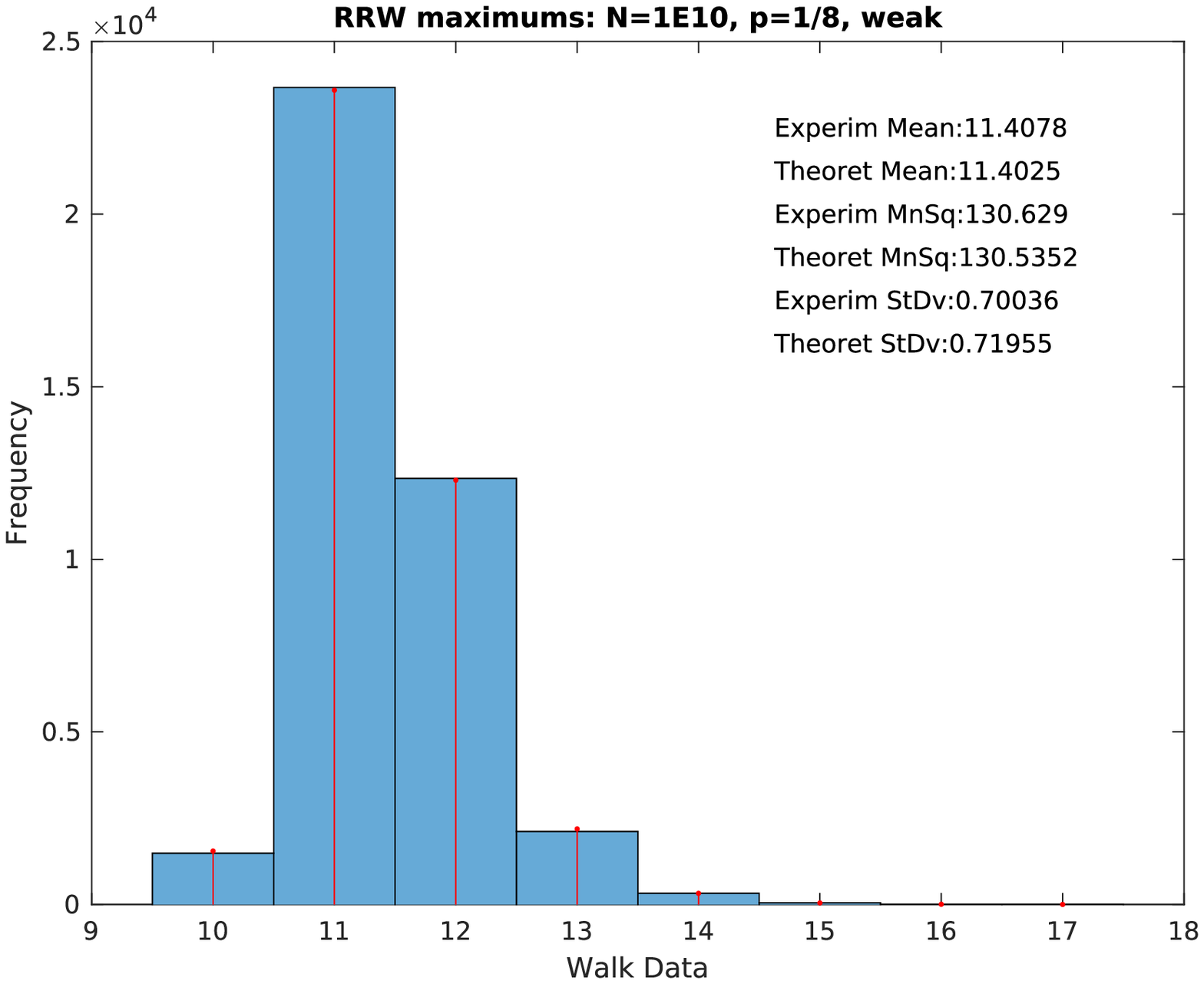}%
\caption{ }%
\end{figure}

Solving these, we have%
\[
\tilde{R}_{k}=\frac{1}{2t}\left[  u\left(  \frac{1-2pqz+t}{2}\right)
^{k}+v\left(  \frac{1-2pqz-t}{2}\right)  ^{k}\right]  ,
\]%
\[
\tilde{Q}_{k}=\frac{1}{2t(1-z)}\left[  u\left(  \frac{1-2pqz+t}{2}\right)
^{k}+v\left(  \frac{1-2pqz-t}{2}\right)  ^{k}-2t\left(  p^{2}z\right)
^{k+1}\right]
\]
where%
\[
u=1+t-\left(  1+3p+t+pt\right)  qz+2pq^{2}z^{2},
\]%
\[
v=-1+t+\left(  1+3p-t-pt\right)  qz-2pq^{2}z^{2}.
\]
Skipping over details, we have
\[
\frac{z_{k}-1}{(p/q)^{2k}}\rightarrow\frac{p(1-2p)^{2}}{q^{3}}%
\]
as $k\rightarrow\infty$ (note the exponent $2k$ in the denominator), which
implies that%
\[
P\left\{  M_{n}\leq\log_{q^{2}/p^{2}}(n)+h\right\}  \sim\exp\left[
-\frac{p(1-2p)^{2}}{2q^{3}}\left(  \frac{q^{2}}{p^{2}}\right)  ^{-h}\right]
\]
as $n\rightarrow\infty$ (note the coefficient $2$ in the denominator).
\ Finally, the discussion in Section 3 applies with $r$ replaced by $r^{2}$.

\section{Queue Data}

Let $n=10^{10}$. \ For each $p\in\{1/5,1/3\}$, we\ generated $40000$ traffic
light queues ($\ell=1$) and produced an empirical histogram for the maximum
$M_{n}$. \ Figures 11--12 contain these histograms (in blue) along with our
theoretical predictions (in red). \ The fit is excellent.

Similarly, we\ generated $40000$ TLQs ($\ell=2$ and $\ell=3$) and produced a
histogram for the maximum $M_{n}$. \ Figures 13--16 contain these histograms.
\ A\ conjecture in \cite{Fi2-mxasr} -- that such distributions do not depend
on the value of $\ell$ -- is evidently false. \ The word \textquotedblleft
theoretical\textquotedblright\ here refers to the ill-informed predictions
emerging from $\ell=1$. \ It would be good someday to understand the true
distributional limits occurring for $\ell\geq2$, even if only approximately.%
\begin{figure}[ptb]%
\centering
\includegraphics[
height=3.0251in,
width=3.691in
]%
{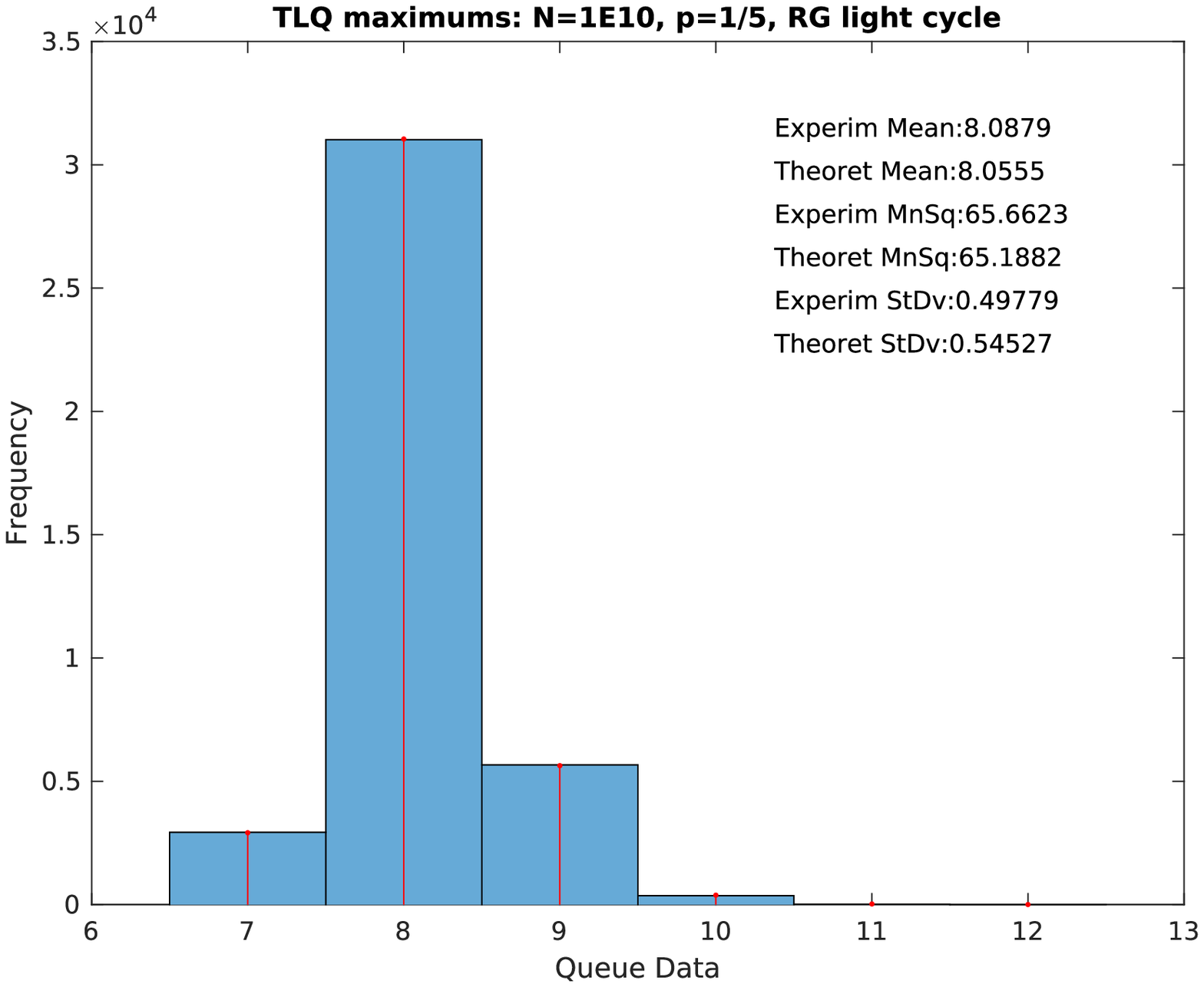}%
\caption{$\ell=1$}%
\end{figure}
\begin{figure}[ptb]%
\centering
\includegraphics[
height=3.0251in,
width=3.691in
]%
{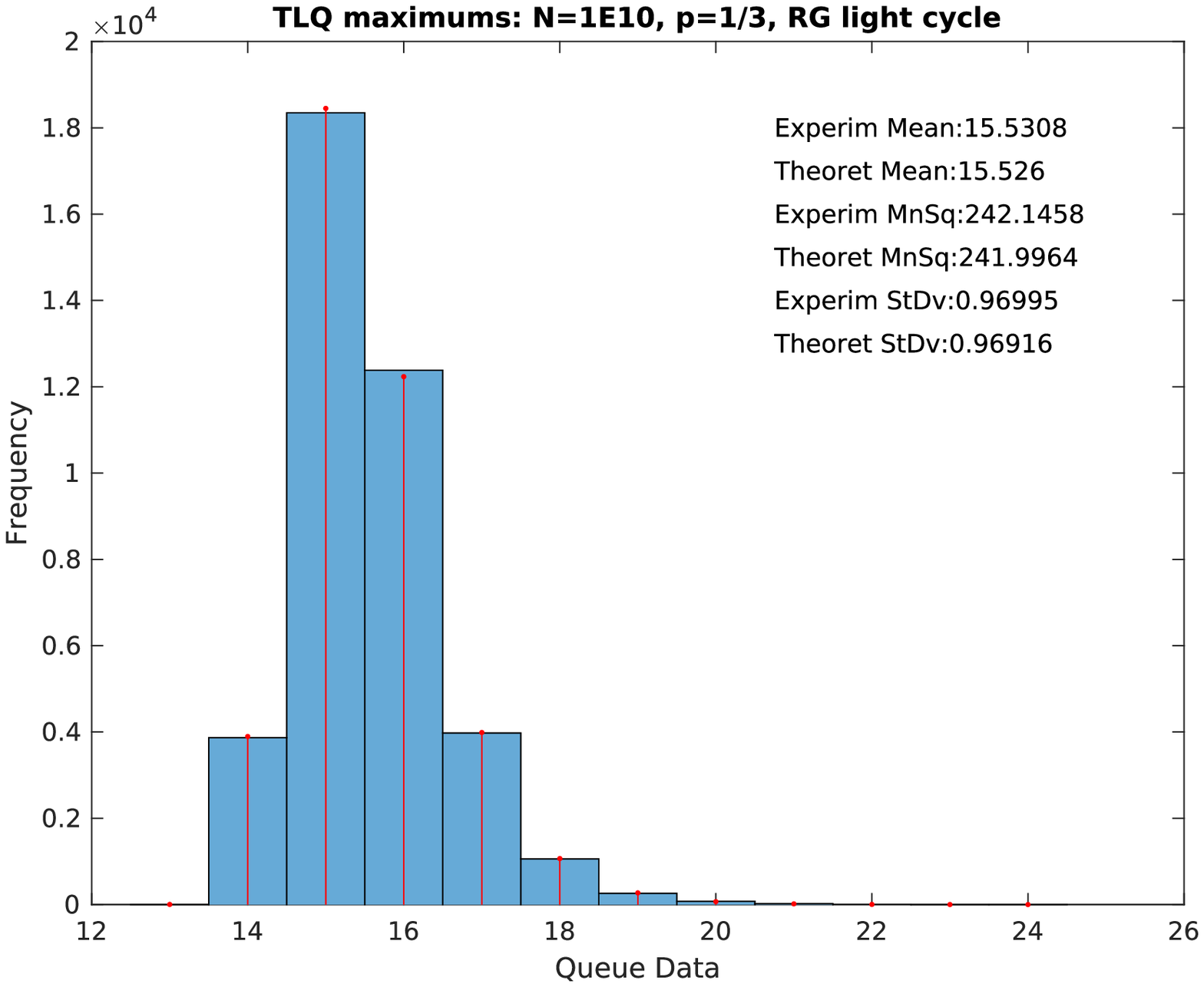}%
\caption{$\ell=1$}%
\end{figure}
\begin{figure}[ptb]%
\centering
\includegraphics[
height=3.0251in,
width=3.691in
]%
{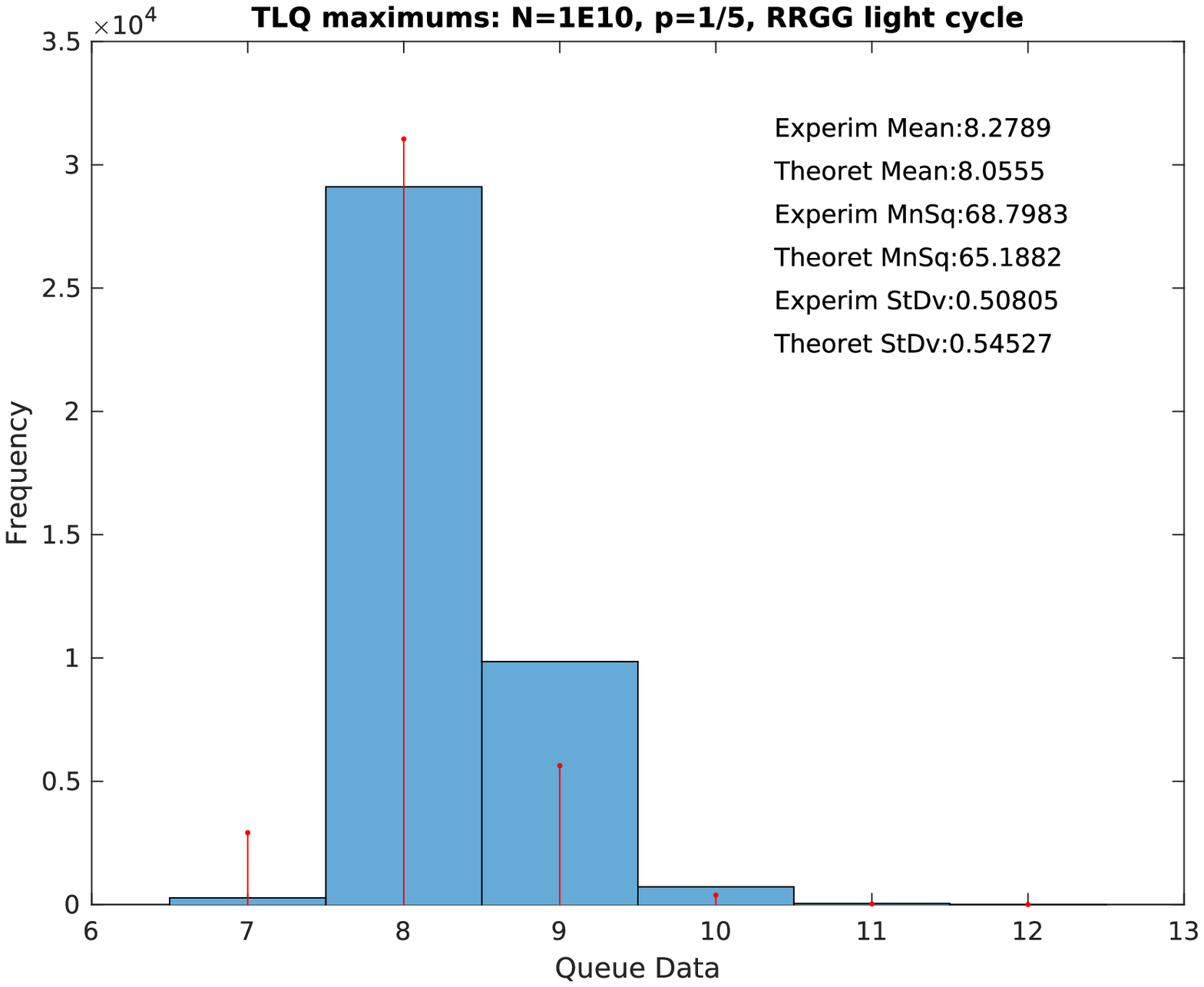}%
\caption{Theory for $\ell=1$ does not carry over to $\ell=2$}%
\end{figure}
\begin{figure}[ptb]%
\centering
\includegraphics[
height=3.0251in,
width=3.691in
]%
{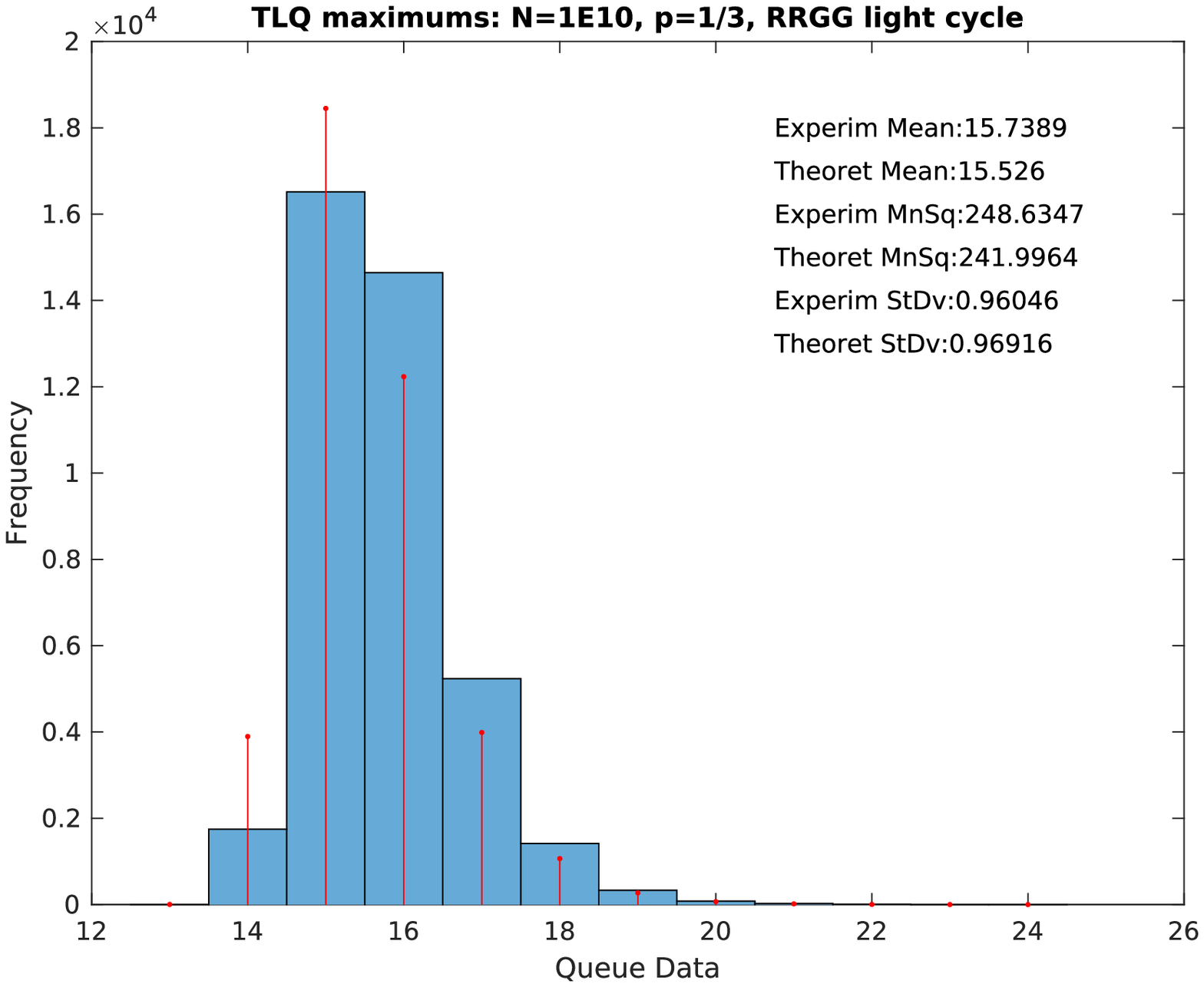}%
\caption{Theory for $\ell=1$ does not carry over to $\ell=2$}%
\end{figure}
\begin{figure}[ptb]%
\centering
\includegraphics[
height=3.0251in,
width=3.691in
]%
{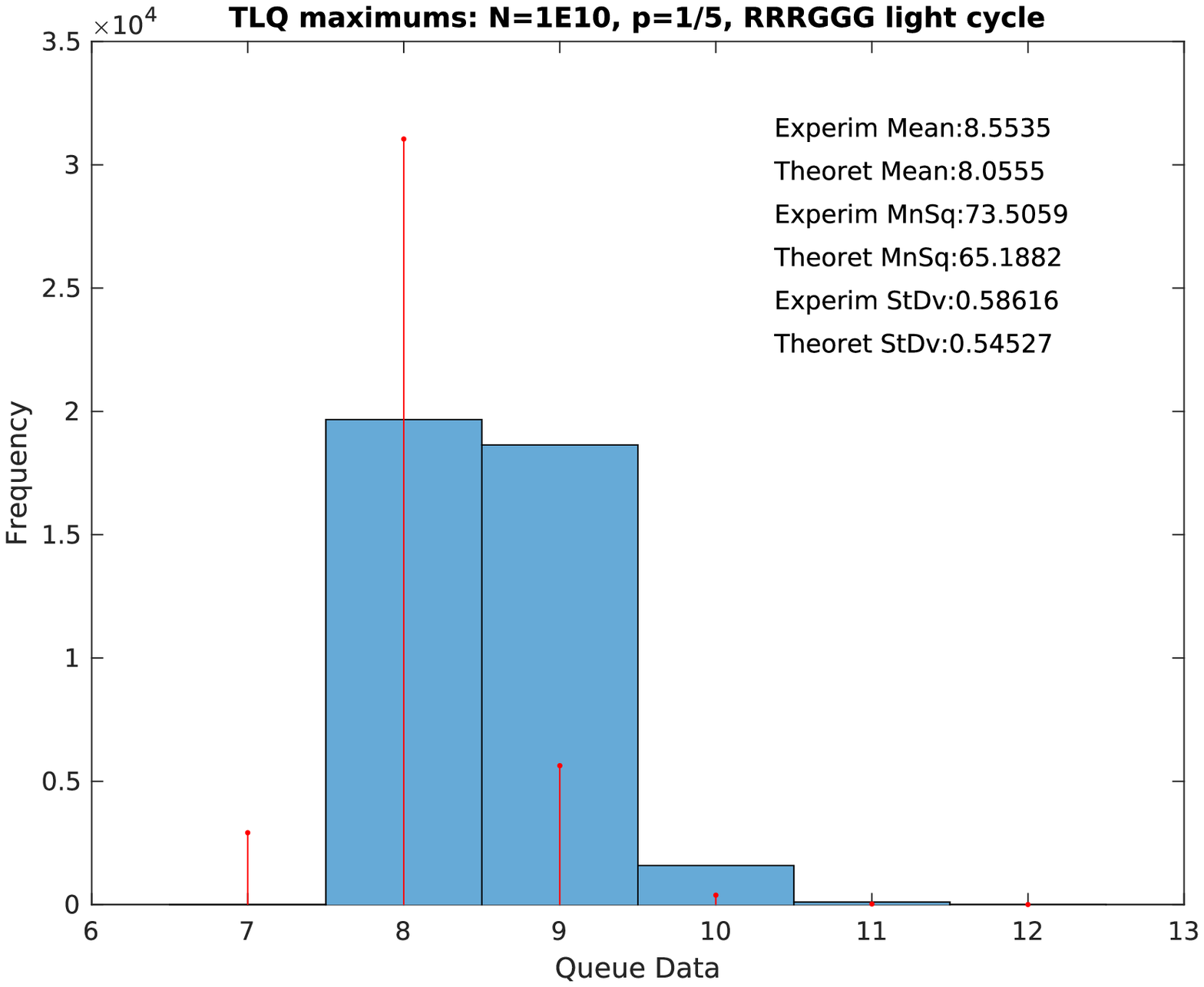}%
\caption{Theory for $\ell=1$ does not carry over to $\ell=3$}%
\end{figure}
\begin{figure}[ptb]%
\centering
\includegraphics[
height=3.0251in,
width=3.691in
]%
{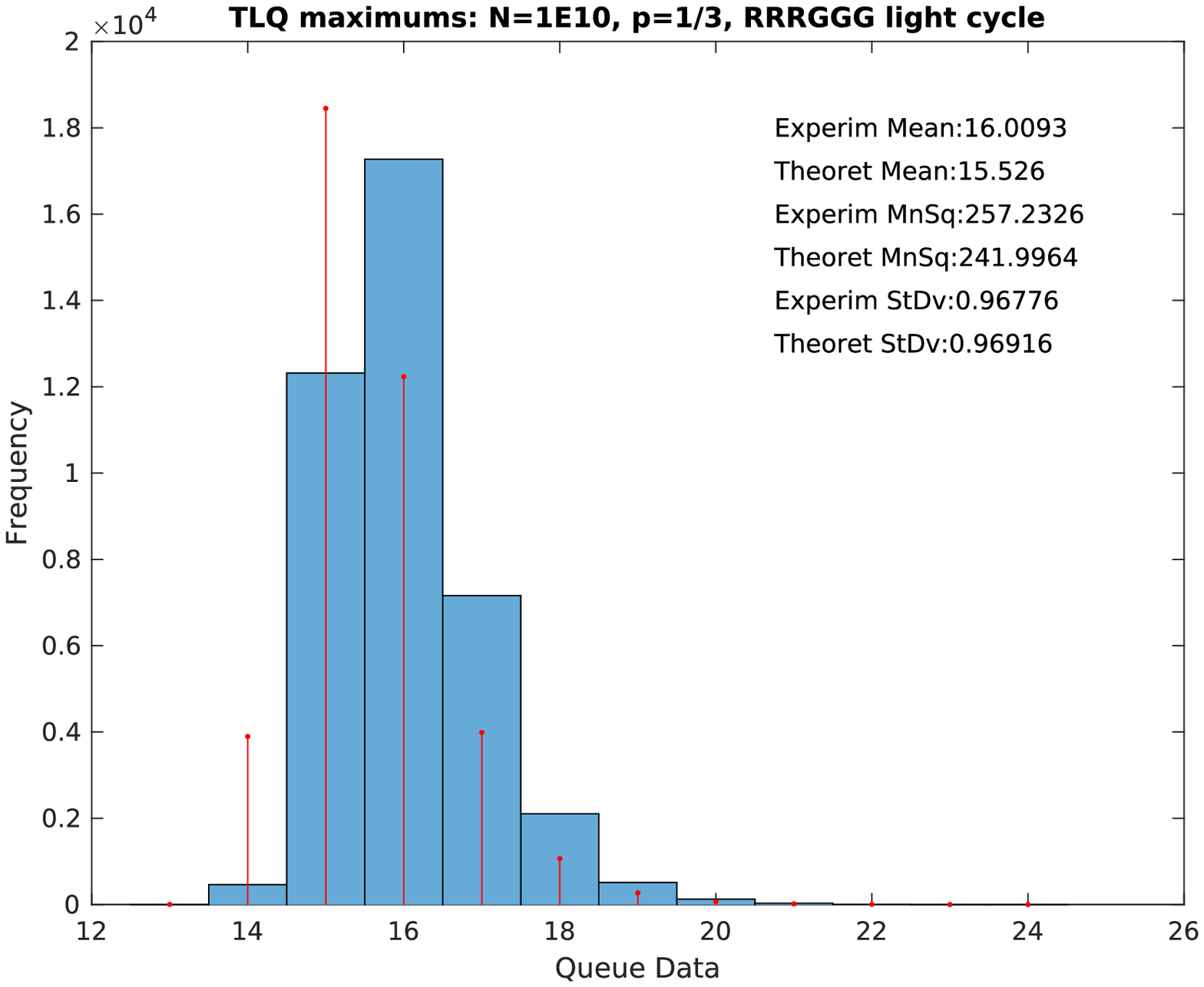}%
\caption{Theory for $\ell=1$ does not carry over to $\ell=3$}%
\end{figure}

\section{Acknowledgements}

I am indebted to Stephan Wagner \cite{HW-mxasr, W-mxasr} for his expertise in
obtaining the discrete Gumbel asymptotics in Section 1. \ Guy Louchard assured
me that his mean/variance formulas \cite{HL-mxasr, LP-mxasr} indeed apply not
just to $q/p=2$ (i.e., $p=1/3$), but to all $q/p>1$ (i.e., $0<p<1/2$); he also
reminded me that, in Section 5, it's best to imagine a $2n$-sequence parsed
into blocks $\{S_{1},S_{2}\},\ldots,\{S_{2n-1},S_{2n}\}$. \ The creators of
Julia, Mathematica and Matlab, as well as administrators of the MIT\ Engaging
Cluster, earn my gratitude every day.


\newpage

\begin{thebibliography}{9}                                                                                                %


\bibitem {Fi1-mxasr}S. R. Finch, How far might we walk at random?, arXiv:1802.04615.

\bibitem {Fi0-mxasr}S. R. Finch, Euler-Mascheroni constant,
\textit{Mathematical Constants}, Cambridge Univ. Press, 2003, pp. 28--40; MR2003519.

\bibitem {Fi2-mxasr}S. R. Finch, Maximum queue length for traffic light with
Bernoulli arrivals, arXiv:1802.04621.

\bibitem {HW-mxasr}H. Prodinger and S. Wagner, Bootstrapping and
double-exponential limit laws, \textit{Discrete Math. Theor. Comput. Sci.} 17
(2015) 123--144; MR3325924.

\bibitem {W-mxasr}S. Wagner, Asymmetric random walk on the integers,
unpublished note (2017).

\bibitem {G-mxasr}Wikipedia contributors, Gumbel distribution,
\textit{Wikipedia, The Free Encyclopedia}, http://en.wikipedia.org/wiki/Gumbel\_distribution.

\bibitem {HL-mxasr}P. Hitczenko and G. Louchard, Distinctness of compositions
of an integer: a probabilistic analysis, \textit{Random Structures Algorithms}
19 (2001) 407--437; MR1871561.

\bibitem {LP-mxasr}G. Louchard and H. Prodinger, Asymptotics of the moments of
extreme-value related distribution functions, \textit{Algorithmica} 46 (2006)
431--467 (long version available at http://www.ulb.ac.be/di/mcs/louchard/); MR2291964.%

\begin{tabular}
[c]{lll}
& Steven Finch & \\
& MIT Sloan School of Management & \\
& Cambridge, MA, USA & \\
& \textit{steven\_finch@harvard.edu} &
\end{tabular}

\end{thebibliography}
\end{document}